\documentclass{article}
\usepackage{mathrsfs,amssymb,amsfonts,amsmath,latexsym,color}
\DeclareMathAlphabet{\mathpzc}{OT1}{pzc}{m}{it}

\usepackage{amsthm}
\usepackage{graphicx}
\usepackage[mathscr]{eucal}
\usepackage[utf8]{inputenc}
\usepackage[colorlinks,
linkcolor=blue,
anchorcolor=green,
citecolor=magenta
]{hyperref}
\usepackage[numbers,sort&compress]{natbib}
\def\msa{\,\mathcal{A}}
\def\msb{\,\mathcal{B}}
\def\msc{\,\mathcal{C}}
\def\msf{\,\mathcal{F}}
\def\msg{\,\mathscr{G}}
\def\msh{\,\mathcal{H}}
\def\msi{\,\mathscr{I}}
\def\msm{\,\mathcal{M}}
\def\msr{\,\mathcal{R}}

\def\mst{\,\mathcal{T}}

\def\mtf{\,\mathbb{F}}

\newcommand{\bn}[2]{\genfrac{(}{)}{0pt}{}{#1}{#2}}
\newcommand{\gn}[2]{\genfrac{[}{]}{0pt}{}{#1}{#2}}

\newtheorem{thm}{Theorem}[section]
\newtheorem{prop}[thm]{Proposition}
\newtheorem{lem}[thm]{Lemma}

\title{Non-trivial intersecting families of finite affine spaces}
\author{Chao Gong\quad Benjian Lv\quad Kaishun Wang\\
	{\footnotesize   \em  Sch. Math. Sci. {\rm \&} Lab. Math. Com. Sys., Beijing Normal University, Beijing, 100875,  China}
}
\date{}
\date{}

\begin{document}

\maketitle
\begin{abstract}
Guo and Xu determined the maximum size of intersecting families over finite affine spaces and showed that any family reaches maximum size must be trivial. In this paper, we characterize non-trivial intersecting family with maximum size. 
\end{abstract}
\section{Introduction}
Let $[n]$ denote the set $\{1,\ldots,n\}$ and, for $0\leq k\leq n$, let $\binom{[n]}{k}$ denote the family of all $k$-subsets of $[n]$. A family $\msf\subseteq\bn{[n]}{k}$ is called \emph{$t$-intersecting} if $|F_{1}\cap F_{2}|\geq t$ for all $F_{1},F_{2}\in\msf$. 
Erd\H{o}s, Ko and Rado \cite{KO} determined the maximum size of a $t$-intersecting family.
\begin{thm}{\em(\cite{KO})}\label{th1.1}
	Let $k\geq t\geq 1$, $n>n_0(k,t)$ and $\msf\subseteq\bn{[n]}{k}$ be a $t$-intersecting family. Then $|\msf|\leq \bn{n-t}{k-t}$. Equality holds iff $\msf$ consists of all $k$-subsets containing a fixed $t$-subset.
\end{thm}

The best value of $n_0(k,t)$ is $(t+1)(k-t+1)$. It was first proved by Frankl \cite{f1978} for $t\geq 15$, and completely determined by Wilson \cite{RMW} based on the eigenvalue method. In \cite{f1978}, Frnakl also gave a conjecture about the maximum size of a $t$-intersecting subfamily of $\bn{[n]}{k}$ for any positive $n,k$ and $t$. This conjecture was partially proved by Frankl and F\"{u}redi \cite{PFZF}, and completed settled by Ahlswede and Khachatrian \cite{RH}. 

A $t$-intersecting family is called \emph{trivial} if every element contains a fixed $t$-subset of $[n]$. Theorem \ref{th1.1} shows that maximum $t$-intersecting subfamily of $\bn{[n]}{k}$ with $n>n_0(k,t)$ must be trivial. Hilton and Milner \cite{AE} determined the maximum size of a non-trivial $1$-intersecting family. The structure of a non-trivial $t$-intersecting subfamily of $\bn{[n]}{k}$ with maximum size was first considered by Frankl \cite{Peter78}, and again, completely determined in \cite{RH}.

Let $V$ denote an $n$-dimensional vector space over the finite field $\mathbb{F}_q$ and $\gn{V}{k}_q$ denote the family of all $k$-dimensional subspaces of $V$. For $n,k\in\mathbb{Z}^+$, define the \emph{Gaussian binomial coefficient} by
$$\gn{n}{k}_q:=\prod_{0\leq i<k}\frac{q^{n-i}-1}{q^{k-i}-1}.$$
Note that the size of $\gn{V}{k}_{q}$ is $\gn{n}{k}_q$. From now on, we will omit the subscript $q$.

A family $\msf\subseteq\gn{V}{k}$ is called \emph{$t$-intersecting} if $\dim(F_1\cap F_2)\geq t$ holds for any $F_1,F_2\in\msf$. A $t$-intersecting family is called \emph{trivial} if every element contains a fixed $t$-space of $V$. The structure of a $t$-intersecting subfamily of $\gn{V}{k}$ with maximum size was partially determined by Hsieh \cite{Hsieh}, Frankl and Wilson \cite{PR}, and completely obtained by Tanaka \cite{TANAKA2006}. Their results showed that maximum $t$-intersecting subfamily of $\gn{V}{k}$ with $\dim(V)\geq 2k+1$ must be trivial.


Suppose $A$ and $B$ are subspaces of $V$. We say $A$ intersects $B$ if $\dim(A\cap B)\geq 1$. Let $\msf\subseteq\gn{V}{k}$ be an $1$-intersecting family. Denote \emph{the covering number $\tau(\msf)$} is the minimum dimension of a subspace of $V$ that intersects all elements of $\msf$. Note that $\msf$ is non-trivial iff $\tau(\msf)\geq 2$. Blokhuis et al. \cite{AB} determined the maximum size of an $1$-intersecting family $\msf$ with $\tau(\msf)\geq 2$.

\begin{thm}{\em{(\cite{AB})}}\label{1.1}
Let $k\geq 3$ and either $q\geq 3$ and $n\geq 2k+1$ or $q=2$ and $n\geq 2k+2$. Let $V$ be an n-dimensional vector space over $\mathbb{F}_q$, then for any intersecting family $\msf\subseteq \gn{V}{k}$ with $\tau(\msf)\geq 2$, we have
$$|\msf|\leq \gn{n-1}{k-1}-q^{k(k-1)}\gn{n-k-1}{k-1}+q^k.$$
Equality holds iff
\begin{itemize}
	\item[\em (i)]$\msf=\{F\in\gn{V}{k}:E\subseteq F,\dim(F\cap U)\geq 1\}\cup\gn{E+U}{k}$ for some $E\in\gn{V}{1}$ and $U\in\gn{V}{k}$ with $E\not\subseteq U$.
	\item[\em (ii)]$\msf=\{F\in\gn{V}{k}:\dim(F\cap S)\geq 2\}$ for some $S\in\gn{V}{3}$ if $k=3$.
\end{itemize} 
\end{thm}

There are results of maximum non-trivial $1$-intersecting subfamily of $\gn{V}{k}$ with $\dim(V)=2k$, see \cite{AAT,Fi}. 

Suppose that $P$ is an $k$-dimensional subspace of $\mathbb{F}_q^n$. A coset of $\mathbb{F}_q^n$ relative to an $k$-dimensional subspace $P$ is called an $k$-flat. The \emph{dimension} of an $k$-flat $U+x$ is defined to be the dimension of the subspace $U$, denoted by $\dim(U+x)$. A flat $F_1$ is said to be \emph{incident} with a flat $F_2$, if $F_1$ contains or is contained in $F_2$. The point set $\mathbb{F}_q^n$ with all the flats and the incidence relation among them defined above is said to be the $n$-dimensional \emph{affine space}, denoted by $AG(n,\mathbb{F}_q)$. Denote by $\msm(k,n)$ the set of all $k$-flats in $AG(n,\mathbb{F}_q)$. Denote by $F_1\cap{F_2}$ the intersection of the flats $F_1$ and $F_2$, and by $F_1\cup{F_2}$ the minimum flat containing both $F_1$ and $F_2$. It is known that the intersection of two flats is either a flat or empty-set. 

A family $\msf\subseteq\msm(k,n)$ is called \emph{$t$-intersecting} if $\dim(F_1\cap F_2)\geq t$ holds for any $F_1,F_2\in\msf$. When $t=1$, we say $\msf$ is intersecting. A $t$-intersecting family is called \emph{trivial} if every element of $\msf$ contains a fixed $t$-flat of $\mathbb{F}_q$. Guo and Xu \cite{GUO2017} determined the maximum size of a $t$-intersecting family and showed that maximum $t$-intersecting family must be trivial.

\begin{thm}{\em(\cite[Theorem 1.3]{GUO2017})}
	Let $n\geq 2k+1$ and $\msf\subseteq\msm(k,n)$ be a $0$-intersecting family. Then $|\msf|\leq\gn{n}{k}$. Equality holds iff $\msf$ consists of all elements containing a fixed vector.
\end{thm}
\begin{thm}{\em(\cite[Theorem 1.4]{GUO2017})}
    Let $t\geq 1$, $n\geq 2k+t+2$ and $\msf\subseteq\msm(k,n)$ be a $t$-intersecting family. Then $|\msf|\leq{\gn{n-t}{k-t}}$. Equality holds iff $\msf$ consists of all elements containing a fixed $t$-flat.
\end{thm}

In 2019, Guo \cite{GUO2019} determined the maximum size of a non-trivial $0$-intersecting family.
\begin{thm}{\em(\cite[Theorem 1.6]{GUO2019})}
	Let $n\geq 2k+1\geq 3$ and $\msf\subseteq\msm(k,n)$ be an non-trivial $0$-intersecting family. Then 
	\[|\msf|\leq 1+\gn{n-1}{k-1}+\sum_{i=0}^k q^{i(i+1)}(q^{k-i}-1)\gn{n-k-1}{i}\gn{k}{1}.\]
	Equality holds iff there exist $U\in\msm(k,n)$ and $x\in \mathbb{F}_q^n\setminus\{x\}$ such that
	\[\msf=\{F\in\msm(k,n):x\in F, F\cap U\neq\emptyset\}\cup\{x\}.\]
\end{thm}

Let $\msf\subseteqq\msm(k,n)$ be an intersecting family. $T\in AG(n,\mathbb{F}_q)$ is called a \emph{covering flat} of $\msf$ if $\dim(T\cap F)\geq 1$ holds for any $F\in\msf$. The covering number $\tau(\msf)$ is the minimum dimension of a covering flat of $\msf$. Note that $\msf$ is trivial iff $\tau(\msf)=1$. In this paper, we determine the maximum size of an intersecting family $\msf\subseteq \msm(k,n)$ with $\tau(\msf)\geq 2$.



For any fixed $n,k\in\mathbb{N}$ and prime power $q$, denote
\[f(n,k,q)=\gn{n-1}{k-1}-q^{k(k-1)}\gn{n-k-1}{k-1}+q^k.\]
Suppose $A,B\in AG(n,\mtf_q)$. Let $A'$ be the subspace related to $A$ and
\[\msh_{A,B}(\gamma)=\{F\in\msm(k,n):A\subseteq F,~\dim(F\cap B)\geq \gamma\}.\]
We say $\msf$ is an HM-type family if there exist $E\in\msm(1,n)$, $U\in\msm(k,n)$ with $\dim(E\cap U)=0$ and $\alpha_1,\ldots,\alpha_{q^k}\in(E\cup U)$ such that
$$\msf=\msh_{E,U}(1)\cup\{A_i+\alpha_i:i\in[q^k]\},$$
where $A_1,\ldots, A_{q^k}$ are the $k$-subspaces of $(E\cup U)'$ not containing $E'$. 

We say $\msf$ is an $\msf_3$-type family if there exist $D\in\msm(3,n)$ with
\[\gn{D'}{2}=\{B_1,\ldots,B_{\gn{3}{2}}\}\]
and $\beta_1,\ldots,\beta_{\gn{3}{2}}\in D$ such that
\[\msf=\bigcup_{R\in \msb}\{F\in\msm(k,n):R\subseteq F\},\]
where
\[\msb=\{B_i+\beta_{i}:i\in[\textstyle\gn{3}{2}]\}.\]
The main result is as follows.

\begin{thm}\label{1.2}
Suppose $k\geq 3$, $n\geq 2k+4$ and $(n,q)\neq(2k+4,2)$. Let $\msf\subseteq \msm(k,n)$ be an intersecting family with $\tau(\msf)\geq 2$. Then $|\msf|\leq f(n,k,q)$. Equality holds iff
\begin{itemize}
\item[\em(i)]$\msf$ is an HM-type family.
\item[\em(ii)]$\msf$ is an $\msf_3$-type family if $k=3$.
\end{itemize}
\end{thm}

In Section $2$, we give some equalities and inequalities which are used to prove Theorem \ref{1.2}. In section $3$, we prove Theorem \ref{1.2}.
\section{Some lemmas}

In this section, we shall give some lemmas which are used to prove Theorem \ref{1.2}. Let $V$ be a space of dimension $(n+l)$ over $\mathbb{F}_{q}$ and $W$ be a fixed $l$-subspace of $V$. A subspace $P\in\gn{V}{m}$ is called an \emph{$(m,k)$-type subspace} of $V$ if $\dim(P\cap W)=k$. Let $N'(m_{1},k_{1};m,k;n+l,n)$ be the number of subspaces of type $(m,k)$ in $V$ containing a given subspace of type $(m_{1},k_{1})$.

\begin{lem}{\emph{(\cite[Lemma 2.3]{KW})}}\label{2.1}
$N^{'}(m_{1},k_{1};m,k;n+l,n)\not= 0$ iff
\[0\leq k_{1}\leq k\leq l,~0\leq m_{1}-k_{1}\leq m-k\leq n.\]
Moreover, if $N'(m_1,k_1;m,k;n+l,n)>0$, then
\[N^{'}(m_{1},k_{1};m,k;n+l,n)=q^{(l-k)(m-k-m_{1}+k_{1})}\gn{n-(m_{1}-k_{1})}{(m-k)-(m_{1}-k_{1})}\gn{l-k_{1}}{k-k_{1}}.\]
\end{lem}

\begin{lem}\label{2.3}{\em(\cite[Theorem 1.18]{WAN})}
The number of $k$-flats in $AG(n,\mtf_q)$ contained in a given $m$-flat, where $0\leq k\leq m\leq n$, is equal to $q^{m-k}\gn{m}{k}$.
\end{lem}
\begin{lem}\label{2.4}{\textup{(\cite[Theorem 1.19]{WAN})}}
The number of $m$-flats in $AG(n,\mtf_q)$ containing a given $k$-flat, where $0\leq k\leq m\leq n$, is equal to $\gn{n-k}{m-k}$.
\end{lem}

\begin{lem}\label{2.5}{\em(\cite[Theorem 1.20]{WAN})}
Let $F_1=F'_1+f_1$ and $F_2=F'_2+f_2$ be two flats of $AG(n,\mathbb{F}_q)$. The following hold.
\begin{itemize}
	\item[\em (i)]$F_1\cap F_2\neq\emptyset$ iff $(f_1-f_2)\in (F'_1+F'_2)$.
	\item[\em (ii)]If $F_1\cap F_2\neq\emptyset$, then $F_1\cap F_2=F'_1\cap F'_2+x$, where $x\in F_1\cap F_2$.
	\item[\em (iii)]$F_1\cup F_2=(F'_1+F'_2)+\langle f_2-f_1\rangle+f_1$. In particular,
	\[
	\dim(F_1\cup F_2)=
	\begin{cases}
	\dim F_1+ \dim F_2-\dim(F_1\cap F_2), &\text{if}~F_1\cap F_2\neq\emptyset,\\
	\dim F_1 + \dim F_2-\dim(F'_1\cap F'_2)+1, &\text{if}~F_1\cap F_2=\emptyset.
	\end{cases}
	\]
\end{itemize}
\end{lem}

For any $A, B\in AG(n,\mathbb{F}_q)$ and $\mathcal{F}\subseteq\mathcal{M}(k,n)$, let
\[\msf'=\{F':F\in\msf\},~\msf_{A}=\{F\in\msf:A\subseteq F\},~\msm_A=\{F\in\msm(k,n):A\subseteq F\}.\]
We say that $A$ intersects $B$ if $\dim(A\cap B)\geq 1$.
\begin{lem}\label{2.6}
Let $\msf\subseteq\msm(k,n)$ be an intersecting family. Then $\msf'$ is an intersecting family of $k$-subspaces with $|\msf'|=|\msf|$ and $\tau(\msf')\leq \tau(\msf)$.
\end{lem}
\begin{proof}
Let $F_1$ and $F_2$ be two distinct flats. By Lemma \ref{2.5}, one gets that $F_1$ intersects $F_2$ only if
\[\dim(F_1'\cap F_2')\geq 1,~F_1'\neq F_2'.\]
Let $T$ be a covering flat of $\msf$ with $\dim(T)=\tau(\msf)$. It follows that $\msf'$ is an intersecting family with $|\msf'|=|\msf|$ and $T'$ is a covering subspaces of $\msf'$, which implies that
\[\tau(\msf')\leq \dim(T')=\dim(T)=\tau(\msf).\] 
\end{proof}

\begin{lem}\label{2.7}
Suppose $m\geq 2$, $M\in\msm(m+1,n)$ and $E\in\msm(1,M)$. Let $A_1,\ldots,A_{\gn{m+1}{m}}$ be the $m$-subspaces of $M'$ and $A_1,\ldots,A_{q^m}$ be those not containing $E'$. 
\begin{itemize}
\item[\em(i)] Let $F_1,F_2$ be two distinct elements of $\msm(m,M)$. Then $F_1$ intersects $F_2$ iff $F_1'\neq F_2'$.
\item[\em(ii)] $\msf\subseteq\msm(m,M)$ is an intersecting family iff  there exist $\alpha_1,\ldots,\alpha_{\gn{m+1}{m}}\in M$ such that
\[\msf\subseteq\{A_i+\alpha_i:i\in[\textstyle\gn{m+1}{m}]\}.\]
\item[\em(iii)] $\msf\subseteq\msm(m,M)$ is an intersecting family with $\dim(F\cap E)=0$ holds for any $F\in\msf$ iff there exist $\beta_1,\ldots,\beta_{q^m}\in M$ such that
\[\msf\subseteq\{A_i+\beta_i:i\in[q^m]\}.\]
\end{itemize}
\end{lem}
\begin{proof}
	It is direct results of Lemma \ref{2.5}.
\end{proof}
\begin{lem}\label{2.8}
Suppose $a<b\leq k$. Let $A\in\msm(a,n)$, $B\in\msm(b,n)$ such that $A$ does not intersect $B$ and $M=A\cup B$.
\begin{itemize}
\item[\em(i)]
We have $b+1\leq \dim(M)\leq a+b+1$ and
\[\dim(M)=a+b+1\Leftrightarrow A\cap B=\emptyset,~A'\cap B'=0.\]
\item[\em(ii)]
One gets that $\msh_{A,B}(1)\subseteq \msh_{A,M}(a+1)$. If $\dim(M)=a+b+1$, then $\msh_{A,B}(1)\subseteq \msh_{A,M}(a+2)$.
\end{itemize}
\end{lem}
\begin{proof}
(i) Note that $0\leq \dim(A'\cap B')\leq a$. The desired results follow by Lemma \ref{2.5}.

\noindent(ii) Let $F\in\msh_{A,B}(1)$. Note that $\dim(F\cap B)\geq 1$ and $A$ does not intersect $(F\cap B)$. One gets that
\begin{equation}\label{eq10}
    A\subsetneqq \left(A\cup (F\cap B)\right)\subseteq (F\cap M),
\end{equation}
which implies that $F\in\msh_{A,M}(a+1)$.

Now assume that $A\cap B=\emptyset$ and $A'\cap B'=0$. Let $F\in\msh_{A,B}(1)$. One gets that $\dim(F\cap B)\geq 1$ and
\[(F\cap B)\cap A=\emptyset,~(F\cap B)'\cap A'=0.\]
By Lemma \ref{2.5} and (\ref{eq10}), it follows that
\[\begin{split}
\dim(F\cap M)\geq \dim\left(A\cup (F\cap B)\right)=\dim(A)+\dim(F\cap B)+1-\dim((F\cap B)'\cap A')\geq a+2,
\end{split}\]
which implies that $F\in\msh_{A,M}(a+2)$.
\end{proof}

\begin{lem}\label{2.12}
Suppose $s<k<n/2$. Let $\mathscr{F}\subseteq \mathscr{M}(k,n)$ be an intersecting family, $U\in\mathscr{F}$ and $S\in\msm(s,n)$ such that $S$ dose not intersect $U$. Denote $M=S\cup U$.
\begin{itemize}
\item[\em{(i)}]
If $\dim(M)\leq k+s$, then there exists an $T\in\msm(s+1,M)$ such that $S\subseteq T$ and $|\msf_S|\leq\gn{k}{1}|\msf_T|$.
\item[\em(ii)]
If $\dim(M)=k+s+1$, then there exists $R\in\msm(s+2,M)$ such that $S\subseteq R$ and $|\msf_S|\leq\gn{k}{1}^2|\msf_R|$.
\item[\em(iii)]
$|\msf_{S}|\leq \gn{k}{1}\gn{n-s-1}{k-s-1}$.
\end{itemize}
\end{lem}
\proof
Denote $m=\dim(M)$ and
\[\begin{split}
\msa=\{A\in\msm(s+1,M):S\subseteq A\},~\msb=\{B\in\msm(s+2,M):S\subseteq B\}.
\end{split}\]
By Lemmas \ref{2.4} and \ref{2.8}, one gets that $\msf_S\subseteq \msh_{S,U}(1)$ and 
\[k+1\leq m\leq k+s+1,~|\msa|=\gn{m-s}{1},~|\msb|=\gn{m-s}{2}.\]

\noindent(i) Note that
\[\msf_S\subseteq \msh_{S,U}(1)\subseteq \msh_{S,M}(s+1)=\bigcup_{A\in\msa}\msm_A\]
by Lemma \ref{2.8}. It follows that
\[\msf_S=\bigcup_{A\in\msa}\msf_A.\]
There exists $T\in\msa$ such that $|\msf_T|$ is maximum. Thus,
\[\begin{split}
|\msf_S|\leq\sum_{A\in\msa}|\msf_A|\leq |\msa||\msf_T|\leq \gn{k}{1}|\msf_T|.
\end{split}\]

\noindent(ii) Similar to the proof of (i), we obtain that
\[\msf_S=\bigcup_{B\in\msb}\msf_B\] 
and there exists $R\in\msb$ such that $|\msf_R|$ is maximum. The fact that
\[|\msb|=\gn{k}{2}<\gn{k}{1}^2\]
implies the desired result.

\noindent(iii) Observe that
\[|\msf_T|\leq \gn{n-s-1}{k-s-1},~|\msf_R|\leq \gn{n-s-2}{k-s-2},~\gn{k}{1}\gn{n-s-2}{k-s-2}<\gn{n-s-1}{k-s-1}\]
by Lemmas \ref{2.3} and \ref{2.10}. It is a direct result of (i) and (ii).
$\qed$

\begin{lem}\label{2.13}
Let $a$, $b$ be positive integers with $a\leq b-2$. Suppose $\{\alpha_n\}$ is a infinite sequence of numbers with $\alpha_n\in\{1,2\}, n\in\mathbb{N}^+$ and $\beta_k=a+\sum_{i=1}^k\alpha_i$, $k\in\mathbb{N}^+$. Then there exists $r$ such that $\beta_{r}\in\{b-1,b\}$.
\end{lem}
\begin{proof}
Let $\msa=\{k\in\mathbb{N}^+: \beta_k>b\}$. Note that $(b-a+1)\in\msa$ and $1\not\in\msa$. By the least number principle, $\msa$ has a minimum element $l$ with $l\geq 2$. Let $r=l-1$. It is routine to check that $\beta_{r}\in\{b-1,b\}$.
\end{proof}

\begin{lem}\label{2.14}
Let $\msf\subseteq \mathscr{M}(k,n)$ be an intersecting family with $\tau(\msf)=t$. Then for any $S\in\msm(s,n)$ with $s\leq t$, 
\[|\msf_{S}|\leq{\gn{k}{1}}^{t-s}\gn{n-t}{k-t}.\]
\end{lem}
\proof
We only need to consider to case when $s<t$. Note that $S$ is not a covering flat of $\msf$. If $s=t-1$, then the desired result follows by Lemma \ref{2.12}. Now assume that $s\leq t-2$. 
By Lemmas \ref{2.12} and \ref{2.13}, there exists a sequence of flats
\[S=S_0\subseteq S_1\subseteq\cdots\subseteq S_r\]
such that $\dim(S_r)\in\{t-1,t\}$ and
\[\dim(S_{i})-\dim(S_{i-1})=\alpha_i\in\{1,2\},~|\msf_{S_{i-1}}|\leq \gn{k}{1}^{\alpha_{i}}|\msf_{S_i}|,~i\in[r].\]
Note that $S_r$ is not a covering flat of $\msf$ if $\dim(S_r)=t-1$. By Lemma \ref{2.12}, one gets that
\[|\msf_S|\leq
\begin{cases}
\gn{k}{1}^{t-s}|\msf_{S_r}|\leq \gn{k}{1}^{t-s}\gn{n-t}{k-t},&\text{if}~\dim(S_r)=t,\\
\gn{k}{1}^{t-s-1}|\msf_{S_r}|\leq \gn{k}{1}^{t-s}\gn{n-t}{k-t},&\text{if}~\dim(S_r)=t-1.
\end{cases}\]
$\qed$
\begin{lem}\label{2.10}{\em(\cite[Lemma 2.1]{AB})}
	Let $a\geq 0$ and $n\geq k\geq a+1$ and $q\geq 2$. Then
	\[\gn{k}{1}\gn{n-a-1}{k-a-1}<\frac{1}{(q-1)q^{n-2k}}\gn{n-a}{k-a}.\]
\end{lem}

\begin{lem}\label{2.11}
	Suppose $k\geq 3$, $q\geq 2$, $r\geq 4$ and $n=2k+r$ with $(r,q)\neq(4,2)$. Then
	\[
	f(n,k,q)>{{k}\brack{1}}{{n-2}\brack{k-2}}-q{{k}\brack{2}}{{n-3}\brack{k-3}}>\left(1-\frac{1}{q^r(q^2-1)}\right)\gn{k}{1}\gn{n-2}{k-2}.
	\]
\end{lem}
\proof
The first inequality is the direct result of Lemma 2.3 in \cite{AB}, and the second inequality follows by
\[q\gn{k}{2}=\frac{(\gn{k}{1}-1)\gn{k}{1}}{q+1}<\frac{\gn{k}{1}^2}{q+1}\]
and Lemma \ref{2.10}.
$\qed$

\section{Proof of Theorem \ref{1.2}}
In this section, we always assume that $k\geq 3$, $n\geq 2k+4$ and $(n,q)\neq (2k+4,2)$. Let $V$ be a $n$-space over $\mathbb{F}_{q}$ and denote $\msm(k,n)$ be the set of all $k$-flats contained in $V$.
\begin{lem}\label{2.9}
	Suppose $E\in\msm(1,n)$ and $U\in\msm(k,n)$ with $\dim(E\cup U)=k+1$. Denote $M=E\cup U$.
	\begin{itemize}
		\item[\em(i)] We have
		\begin{equation}\label{k1}
			|\msh_{E,M}(2)|=\gn{n-1}{k-1}-q^{k(k-1)}\gn{n-k-1}{k-1}.
		\end{equation} 
		\item[\em(ii)] If $E\cap U\neq\emptyset$, then $\msh_{E,U}(1)=\msh_{E,M}(2)$.
		\item[\em(iii)] If $E\cap U=\emptyset$, then
		\[|\msh_{E,U}(1)|\leq |\msh_{E,M}(2)|-q^{(k-1)(k-2)}\gn{n-k-1}{k-2}.\]
	\end{itemize}
\end{lem}
\proof
\noindent(i) It is routine to check that $\msh_{E,M}(2)$ is an intersecting family and
\[|\msh'_{E,M}(2)|=\gn{n-1}{k-1}-N'(1,1;k,1;n,n-k-1)\]
by Lemma \ref{2.1}, which implies that (\ref{k1}) holds by Lemma \ref{2.6}.

\noindent(ii) Note that $\dim(M')=k+1$ and $U\in\gn{M'}{k}$. It is routine to check that
\[F\cap U\neq\emptyset,~\dim(F'\cap M')=\dim(F\cap M)\geq 2\]
hold for any $F\in\msh_{E,M}(2)$. One gets that
\[\dim(F\cap U)=\dim(F'\cap U')=\dim((F'\cap M')\cap U')\geq 1\]
by Lemma \ref{2.6}, which implies that $F\in\msh_{E,U}(1)$. Together with $\msh_{E,U}(1)\subseteq\msh_{E,M}(2)$ by Lemma \ref{2.8}, the desired result follows.

\noindent(iii) Assume that $E=E'+e$ and $U=U'+u$. One gets that 
\[E'\subseteq U'\subseteq M',~(e-u)\not\in U',~M=M'+e.\]
Pick $A\in\gn{U'}{2}$ with $E'\subseteq A$ and let $S=A+e$. It follows that
\[E\subseteq S\subseteq M,~S\cap U=\emptyset.\]
Denote
\[\msa=\{F\in\msm(k,n): F\cap M=S\}.\]
It is routine to check that
\[\msa\subseteq \msh_{E,M}(2)\setminus\msh_{E,U}(1).\]
Note that $\msa$ is an intersecting family. It follows that
\[\begin{split}
|\msa|&=|\msa'|=\left|\left\{K\in\gn{V}{k}:K\cap M'=A\right\}\right|\\
&=N'(2,2;k,2;n,n-(k+1))=q^{(k-1)(k-2)}\gn{n-k-1}{k-2}.
\end{split}\]
by Lemmas \ref{2.1} and Lemma \ref{2.5}, which implies the desired result by (\ref{k1}).


\begin{prop}\label{3.1}
	Let $\msf\subseteq \msm(k,n)$ be an HM-type family. Then $\msf$ is an intersecting family with $\tau(\msf)=2$ and $|\msf|=f(n,k,q)$.
\end{prop}
\begin{proof}
There exist $E\in\msm(1,n)$, $U\in\msm(k,n)$ with $\dim(E\cap U)=0$ and $\alpha_1,\ldots,\alpha_{q^k}\in(E\cup U)$ such that
\[\msf=\msh_{E,U}(1)\cup\{A_i+\alpha_i:i\in[q^k]\},\]
where $A_1,\ldots,A_{q^k}$ are the $k$-subspaces of $(E\cup U)'$ not containing $E'$. Denote $M=E\cup U$ and
\[\msa=\{A_i+\alpha_i:i\in[q^k]\},~\msb=\{S\in\msm(2,M):E\subseteq S\}.\]
By Lemmas \ref{2.1}, \ref{2.7} and \ref{2.9}, it is routine to check the following hold:
\begin{itemize}
	\item[(i)] For any $F\in\msh_{E,U}(1)$, there exists $S\in\msb$ such that $S\subseteq F$.
	\item[(ii)] Every element of $\msb$ is a covering flat of $\msh_{E,U}(1)$ and $\msa$.
	\item[(iii)] $\msh_{E,U}(1)$ and $\msa$ are intersecting families of $k$-flats with $|\msf|=f(n,k,q)$ and $\tau(\msa')=2$.	
\end{itemize}
It follows that $\msf$ is an intersecting family and every element of $\msb$ is a covering flat of $\msf$, which implies that
\[2\geq \tau(\msf)\geq \tau(\msa)\geq \tau(\msa')=2\]
by Lemma \ref{2.6}
\end{proof}

\begin{prop}\label{3.2}
Let $\msf\subseteq \msm(3,n)$ be an $\msf_3$-type family. Then $\msf$ is an intersecting family with $\tau(\msf)=2$ and $|\msf|=f(n,3,q)$.
\end{prop}
\proof
There exist $D\in\msm(3,n)$ and $\beta_1,\ldots,\beta_{\gn{3}{2}}\in D$ such that $\msf=\bigcup_{A\in\msa}\msm_A$ with
\[\gn{D'}{2}=\{B_1,\ldots,B_{\gn{3}{2}}\},~\msa=\{B_1+\beta_i:i\in[\textstyle\gn{3}{2}]\}.\]
Note that $\msa$ is an intersecting family by Lemma \ref{2.7}. It follows that $\msf$ is an intersecting family and every element of $\msa$ is a covering flat of $\msf$. One gets that 
\[\msf'=\left\{K\in\gn{V}{3}: \dim(K\cap D')\geq 2\right\}\]
which implies that $|\msf|=|\msf'|=f(n,3,q)$ and
\[2\geq \tau(\msf)\geq \tau(\msf')=2\]
by Theorem \ref{1.2} and Lemma \ref{2.6}.
$\qed$

For $\msf\subseteq \msm(k,n)$ with $\tau(\msf)=t$, denote $\mst_{\msf}$ be the family of all $t$-covering flats of $\msf$. We will omit the subscript $\msf$.
\begin{lem}\label{3.3}
Let $\msf\subseteq \msm(k,n)$ be an intersecting family with $\tau(\msf)=2$.
\begin{itemize}
    \item[\em(i)] If there exist $A,B\in\mst$ such that $A$ does not intersect $B$, then $\msf$ is not maximal.
    \item[\em(ii)] If $\msf$ is maximal, then
    \[\bigcup_{T\in\mst}\msm_T\subseteq\msf.\]
\end{itemize}
\end{lem}
\proof
(i) Let $A=A'+a$ and $B=B'+b$. If $A'\cap B'=0$, then there exists $K\in\gn{V}{k}$ such that
\[A'\subseteq K,~K\cap B'=0,\]
which implies that $(K+a)$ is contained in $\msm_A$ and does not intersect $B$. If $\dim(A'\cap B')\geq 1$, then $A\cap B=\emptyset$. Denote $E=\langle a-b\rangle$. It follows that
\[E\cap (A'+B')=0,~\dim(A'+B')\leq 3.\]
By Lemma \ref{2.1}, there exists $X\in\gn{V}{n-1}$ such that $(A'+B')\subseteq X$ and $X\cap E=0$. Pick $K\in\gn{X}{k}$ with $A'\subseteq K$. The fact that $(B'+K)\cap E=0$ implies that 
\[(K+a)\in\msm_A,~(K+a)\cap B=\emptyset.\]

Let $\msf_1=\msf\cup\msm_A$. It is obvious that $\msf_1$ is an intersecting family and $B$ is not a covering flat of $\msf_1$, which implies that $\msf\subsetneqq\msf_1$. 

\noindent(ii) By (i), one gets that $|\mst|=1$ or $\mst$ is an intersecting family. Let $\msf_2=\msf\cup\bigcup_{T\in\mst}\msm_T$. It is routine to check that $\msf_2\subseteq\msm(k,n)$ is an intersecting family with $\tau(\msf_2)=\tau(\msf)$, which implies that $\msf_2=\msf$ by the maximality of $\msf$.
$\qed$

By Lemma \ref{3.3}, maximal intersecting families with covering number no less than $2$ can be divided into four cases:
\begin{itemize}
\item $\tau(\msf)=2$ and $|\mst|=1$.

\item $\tau(\msf)=2$ and $\mst$ is an intersecting family with $\tau(\mst)=1$.

\item $\tau(\msf)=2$ and $\mst$ is an intersecting family with $\tau(\mst)=2$.

\item $\tau(\msf)\geq 3$.
\end{itemize}
\begin{prop}\label{3.4}
Let $\msf\subseteq\msm(k,n)$ be an intersecting family with $\tau(\msf)=2$ and $\mst=\{T\}$. Then
\begin{equation*}
|\msf|\leq \gn{n-2}{k-2}+q(q+1)\left(\gn{k}{1}-1\right)\gn{k}{1}\gn{n-3}{k-3}.
\end{equation*}

\end{prop}

\begin{proof}
	Denote $\msm(1,T)=\{E_1,\ldots,E_{q(q+1)}\}$. Note that $\tau(\msf)=2$. Then for any $i\in[q(q+1)]$, there exists $U_i\in\msf$ such that $E_i$ does not intersect $U_i$. Let
	\[\begin{split}
	M_i=(E_i\cup U_i),~\msa_i=\{S\in\msm(2,M_i):E_i\subseteq S_i\}\setminus\{T\},~\msb_i=\{R\in\msm(3,M_i):E_i\subseteq R,T\not\subseteq R\}.
	\end{split}\] 
	It is routine to check that
	\[\begin{split}
	\dim(M_i)\in\{k+1,k+2\},~T=E_i\cup(U_i\cap T)\subseteq E_i\cup U_i
	\end{split}\]
	hold for any $i\in[q(q+1)]$. We claim that
	\begin{align}
		\msf_{E_i}\setminus\msf_T\subseteq
		\begin{cases}
		\displaystyle\bigcup_{S\in\msa_i}\msf_S&\text{if} \dim(M_i)=k+1,\\
		\displaystyle\bigcup_{R\in\msb_i}\msf_R&\text{if} \dim(M_i)=k+2
		\end{cases}
	\end{align}
	hold for any $i\in[q(q+1)]$. If $\dim(M_i)=k+1$, then by Lemma \ref{2.8}, one gets that
	\[\left(\msf_{E_i}\setminus\msf_T\right)\subseteq \msf_{E_i}\subseteq\msh_{E_i,U_i}(1)\subseteq \msh_{E_i,M_i}(2).\]
	Therefore, for any $F\in\left(\msf_{E_i}\setminus\msf_T\right)$, there exists $2$-flat $S$ such that $E_i\subseteq S\subseteq M_i$ and $S\subseteq F$. Note that $T\not\subseteq F$. It follows that $S\in\msa_i$, which implies that 
\[\msf_{E_i}\setminus\msf_T\subseteq\bigcup_{S\in\msa_i}\msf_S.\]
If $\dim(M_j)=k+2$, then by Lemma \ref{2.8}, one gets that
\[\left(\msf_{E_j}\setminus\msf_T\right)\subseteq \msf_{E_j}\subseteq \msh_{E_j,U_j}(1)\subseteq \msh_{E_j,M_j}(3).\]
Therefore, for any $F\in\left(\msf_{E_j}\setminus\msf_T\right)$, there exists $3$-flat $E_j\subseteq R\subseteq M_j$ such that $R\subseteq F$. Note that $T\not\subseteq F$. It follows that $R\in\msb_j$, which implies that
\[\msf_{E_j}\setminus\msf_T\subseteq\bigcup_{R\in\msb_j}\msf_R.\]

It obvious that
\[\msf=\bigcup_{i=1}^{q(q+1)}\msf_{E_i}=\bigcup_{i=1}^{q(q+1)}\left(\msf_{E_i}\setminus\msf_T\right)\cup\msf_T\]
and $\msa_i\cap \mst=\emptyset$ for any $i\in[q(q+1)]$. By Lemmas \ref{2.5} and \ref{2.10}, one gets that 
\[|\msa_i|=\gn{k}{1}-1,~|\msf_S|\leq \gn{k}{1}\gn{n-3}{k-3}\]
for any $S\in\msa_i$ if $\dim(M_i)=k+1$, and
\[|\msb_j|=\gn{k+1}{2}-\gn{k}{1},~|\msf_R|\leq \gn{n-3}{k-3}\]
for any $R\in\msb_j$ if $\dim(M_j)=k+2$. The fact that
\[|\msf_T|\leq \gn{n-2}{k-2},~\gn{k+1}{2}-\gn{k}{1}<\gn{k}{1}\left(\gn{k}{1}-1\right)\]
imply the desired result.
\end{proof}

%
\begin{prop}\label{3.5}
Let $\msf\subseteq\msm(k,n)$ be an intersecting family with $\tau(\msf)=2$ and $\mst$ is an intersecting family with $|\mst|\geq 2$ and $\tau(\mst)=1$. Denote
\[m=\dim\left(\bigcup_{T\in\mst}T\right)-1.\]
Then $2\leq m\leq k$ and the following hold.
\begin{itemize}
	\item[\em(i)]If $3\leq m<k$, then
	\[|\msf|\leq\gn{m}{1}\gn{n-2}{k-2}+\left(\gn{k}{1}-\gn{m}{1}\right)\gn{k}{1}\gn{n-3}{k-3}+(q^{m+1}+q^m-1)\gn{n-m}{k-m};\]
    if $m=2$, then
	\[|\msf|\leq(q+1)\gn{n-2}{k-2}+\left(\gn{k}{1}-q-1\right)\gn{k}{1}\gn{n-3}{k-3}+(q^3+q^2-1)\gn{k}{1}\gn{n-3}{k-3};\]
	\item[\em(ii)]If $m=k$, then $|\msf|\leq f(n,k,q)$ and equality holds iff $\msf$ is an HM-type family. 
\end{itemize}

\end{prop}
\proof
Let $\mst=\{T_1,\ldots,T_l\}$ with $l\geq 2$. The fact that 
\[m+1=\dim(M)\geq \dim(T_1\cup T_2)=3\]
implies that $m\geq 2$. Pick $U\in\msf\setminus\msf_E$. Denote $E_i=U\cap T_i$ for any $i\in[l]$ and $X=\bigcup_{i=1}E_i$. It is routine to check that
\begin{equation}\label{eq1}
	E_i\in\msm(1,n),~T_i=E\cup E_i,~E_i'\subseteq X'\subseteq U'
\end{equation}
hold for any $i\in[l]$ and
\[X\subseteq (U\cap M)\subsetneqq M,~M=E\cup X\subseteq E\cup U,\]
which implies that
\begin{equation}\label{eq2}
	\begin{split}
	\dim(X)\leq \dim(U\cap M)<\dim(M)\leq \dim(E\cup U).
	\end{split}
\end{equation}
If $E\cap U\neq\emptyset$, then by Lemma \ref{2.6}, (\ref{eq1}) and (\ref{eq2}), one gets that
\begin{equation}\label{eq3}
	E'\cap E_1'=E'\cap X'=E'\cap U'=0,~E\cap E_1\neq \emptyset,~E\cap X\neq\emptyset.
\end{equation}
If $E\cap U=\emptyset$, then by Lemma \ref{2.6}, (\ref{eq1}) and (\ref{eq2}), one gets that
\begin{equation}\label{eq4}
E\cap E_1=E\cap X=\emptyset,~E'=E_1'\subseteq X'\subseteq U'.
\end{equation}
By Lemma \ref{2.6}, Both (\ref{eq3}) and (\ref{eq4}) imply that
\begin{equation}\label{k0}
	\dim(E\cup U)=k+1,~m\leq k,~U\cap M=X\in\gn{M}{m}.
\end{equation} 

\noindent(i) Pick a fixed $U\in\msf\setminus\msf_E$. Denote $N=E\cup U$ and
\[\begin{split}
\msa&=\{A\in\msm(2,M):E\subseteq A\},\\
\msb&=\{B\in\msm(2,N):E\subseteq B\not\subseteq M\},\\
\msc&=\{C\in\msm(m,M):E\not\subseteq C\}.
\end{split}\]
By (\ref{k0}), Lemmas \ref{2.3} and \ref{2.4}, one gets that
\begin{equation}\label{eq5}
\begin{split}
    &\msf_E=\bigcup_{A\in\msa}\msf_A\cup\bigcup_{B\in\msb}\msf_B,~\msf\setminus\msf_E=\bigcup_{C\in\msc}\msf_C,\\
    &|\msa|=\gn{m}{1},~|\msb|=\gn{k}{1}-\gn{m}{1},~|\msc|=q^{m+1}+q^m-1.
\end{split}
\end{equation}
Note that 
\[\mst\cap\msb=\mst\cap\msc=\emptyset.\]
The desired upper bounds follow by (\ref{eq5}), Lemmas \ref{2.4} and \ref{2.12}.

\noindent(ii) By (\ref{k0}), one gets that $\msf\setminus\msf_E\subseteq\msm(k,M)$. Let $A_1,\ldots,A_{q^k}$ be the elements of $\gn{M'}{k}$ not containing $E'$. If $U\cap E\neq\emptyset$ for any $U\in\msf\setminus\msf_E$, then there exist $\alpha_1,\ldots,\alpha_{q^k}\in M$ such that
\[\msf\setminus\msf_E\subseteq \{A_i+\alpha_i:i\in[q^k]\}\]
by Lemma \ref{2.7}. Pick a fixed $U\in\msf\setminus\msf_E$. One gets that
\[M=E\cup U,~\msf_E\subseteq \msh_{E,U}(1),\]
which implies that $\msf$ is an HM-type family by the maximality of $\msf$ and Proposition \ref{3.1}.

Now assume that there exists an $U\in\msf\setminus\msf_E$ such that $E\cap U=\emptyset$. Observe that $k\geq 3$, $n\geq 2k+4$ and $\left(\msf\setminus\msf_E\right)\subseteq \msm(m,(E\cup U))$ is an intersecting family. One gets that
\[\begin{split}
|\msf_E|&\leq |\msh_{E,U}|\leq f(n,k,q)-q^k-q^{(k-1)(k-2)}\gn{n-k-1}{k-2}<f(n,k,q)-2q^k,\\
|\msf\setminus\msf_E|&\leq\gn{k+1}{1}=q^k+\gn{k}{1}<2q^k
\end{split}\]
by Lemmas \ref{2.7} and \ref{2.9}, which implies that $|\msf|<f(n,k,q)$.
$\qed$
\begin{prop}\label{3.6}
Let $\msf\subseteq\msm(k,n)$ be an intersecting family with $\tau(\msf)=2$ and $\mst$ be an intersecting family with $\tau(\mst)=2$. 
\begin{itemize}
	\item[\em(i)] There exists $D\in\msm(3,n)$ such that $\mst\subseteq \msm(2,D)$ and
	\[\dim(F\cap D)\geq 2,~F\in\msf.\]
	\item[\em(ii)] One gets that
	\begin{equation}\label{eq6}
	|\msf|\leq \gn{3}{2}\left(\gn{n-2}{k-2}-\gn{n-3}{k-3}\right)+\gn{n-3}{k-3}.
	\end{equation}
\end{itemize}
\end{prop}
\proof
Pick $A,B\in\mst$. Let $E=A\cap B$ and $D=A\cup B$. One gets that $\dim(E)=1$, $\dim(D)=3$ and there exists $C\in\mst$ such that $E$ dose not intersect $C$. Denote $E_1=C\cap A$ and $E_2=C\cap B$. It is routine to check that
\[C=E_1\cup E_2,~A=E\cup E_1,~B=E\cup E_2,~D=E\cup C.\] 

\noindent(i) It suffices to prove that if a flat $U$ intersects $A$, $B$ and $C$, then $\dim(U\cap D)\geq 2$. Assume that $E\subseteq U$. Note that $\dim(U\cap C)\geq 1$ and $E$ dose not intersect $(U\cap C)$. It follows that
\[E\subsetneqq E\cup (U\cap C)\subseteq (U\cap D),\]
which implies that $\dim(U\cap D)\geq 2$. Assume that $E\not\subseteq U$. Let $A_1=U\cap A$ and $B_1=U\cap B$. It is routine to check that
\[A_1,B_1\in\msm(1,D),~A_1\neq B_1,~(A_1\cup B_1)\subseteq U\cap D,\]
which implies that 
\[\dim(U\cap D)\geq \dim(A_1\cup B_1)=2.\]

\noindent(ii) Note that $\msf$ is an intersecting family. One gets that $|\msf|=|\msf'|$ and
\[\msf'\subseteq \left\{K\in\gn{V}{k}: \dim(K\cap D)\geq 2\right\}\]
by Lemma \ref{2.6} and (i), which implies that (\ref{eq6}) holds.
$\qed$ 

\begin{prop}\label{3.7}
	Let $\msf\subseteq\msm(3,n)$ be an intersecting family with $\tau(\msf)=2$ and $\mst$ be an intersecting family with $\tau(\mst)=2$. Then $|\msf|\leq f(n,3,q)$ and equality holds iff $\msf$ is an $\msf_3$-type family.
\end{prop}
\proof
By Proposition \ref{3.6}, there exists $D\in\msm(3,n)$ such that $\mst\subseteq\msm(2,D)$ and $\dim(F\cap D)\geq 2$ holds for any $F\in\msf$. Let $A_1,\ldots,A_{\gn{3}{2}}$ be the $2$-subspaces of $D'$. 

If $F\cap D\in\mst$ holds for any $F\in\msf\setminus\{D\}$, then there exists $\alpha_1,\ldots,\alpha_{\gn{3}{2}}\in D$ such that
\[\mst\subseteq\{A_i+\alpha_i:i\in[\textstyle\gn{3}{2}]\}\]
by Lemma \ref{2.7}. Note that
\[\msf\subseteq\bigcup_{T\in\mst}\msf_T\cup\{D\}\subseteq\bigcup_{T\in\mst}\msm_T.\]
One gets that $\msf$ is an $\msf_3$-type family by the maximality of $\msf$.

Assume that there exists $U\in\msf$ such that $U\neq D$ and $U\cap D\not\in\mst$. Let
\[
\msf_1=\{F\in\msf: F\neq D,~(F\cap D)\not\in \mst\},~\mst_1=\{F\cap D:F\in\msf_1\}
\]
and for any $S\in\mst'_1$, denote
\[\msa_S=\{F\in\msf_1: (F\cap D)'=S\}.\]

First, we prove that
\begin{equation}\label{k2}
	\mst'\cap\mst'_1=\emptyset
\end{equation}
by contradiction. Assume that $T\in\mst'\cap\mst_1'$. Then there exist $T_1,T_2\in\msm(2,D)$ and $U\in\msf_1$ such that
\[T_1'=T_2'=T,~U\cap D=T_1,~T_2\in\mst.\]
The fact that $T_1\cap T_2=\emptyset$ implies that $U\cap T_2=\emptyset$, which leads to a contradiction. 

Next, we claim that for any $S\in\mst'_1$, there exits $X_S\in\msm(4,n)$ such that
\begin{equation}\label{k3}
S\subseteq X_S',~\msa_S\subseteq \msm(3,X_S).
\end{equation}
Pick $S\in\mst'_1$. There exist $S_1\in\mst_1$ and $F_1\in\msf_1$ such that 
\[F_1\cap S=S_1,~S_1'=S.\]
Since $S_1\in\mst_1$ is not a covering flat of $\msf$, there exists $F_2\in\msf$ such that $S_1$ dose not intersect $F_2$. Let $S_2=F_2\cap D$. It is routine to check that
\[F_2\neq D,~S_2\in\mst_1,\]
which implies that $F_2\in\msf_1$. Note that $S_1$ dose not intersect $S_2$. One gets that
\[S_2'=S_1'=S\]
by Lemma \ref{2.6}, which implies that
\[F_1'\cap F_2'=S.\]
Let $S_0=F_1\cap F_2$ and $X_S=F_1\cup F_2$. It is routine to check that $\dim(X_S)=4$, $S_0'=S$, $S'\subseteq X_S'$ and
\[\begin{split}
    &S_1\cap S_2=S_1\cap S_0=S_2\cap S_0=\emptyset,\\
	&D=S_1\cup S_2,~F_1=S_0\cup S_1,~F_2=S_0\cup S_2,\\
	&X_S=S_0\cup S_1\cup S_2=S_0\cup D=S_1\cup F_2=S_2\cup F_1,\\
	&S_0\cap D=S_1\cap F_2=S_2\cap F_1=\emptyset.
\end{split}\]
Note that $\{F_1,F_2\}\subseteq\msm(3,X_S)$. Let $U\in\msa_S\setminus\{F_1,F_2\}$ and $S_3=U\cap D$. If $S_0\subseteq U$, then the fact that $S_0\cap S_3=\emptyset$ implies that
\[U=S_0\cup S_3\subseteq S_0\cup D=X_S.\]
If $S_0\not\subseteq U$, then $S_0\cap U=\emptyset$. Let $U_1=F_1\cap U$ and $U_2=F_2\cap U$. It is routine to check that
\[U_1'=U_2'=S,~U_1\neq U_2,\]
which implies that $U_1\cap U_2=\emptyset$ and
\[U=U_1\cup U_2\subseteq F_1\cup F_2=X_S.\]

Note that $|\mst|=|\mst'|$ by Lemma \ref{2.6}. It is routine to check that $\mst'\cup\mst'_1\subseteq \gn{D'}{2}$ and
\begin{equation}\label{eq9}
\msf=\msf_1\cup\{D\}\cup\left(\bigcup_{T\in\mst}\msf_T\setminus\{D\}\right)=\bigcup_{S\in\mst'_1}\msa_S\cup \left(\bigcup_{T\in\mst}\msf_T\setminus\{D\}\right)\cup \{D\}.
\end{equation}
One gets that $|\mst|+|\mst_1'|\leq \gn{3}{2}$ and
\[\msa'_S\subseteq \left\{K\in\gn{X_S}{3}: S\subseteq K\right\},~|\msf_T\setminus\{D\}|\leq \gn{n-2}{1}-1\]
hold for any $T\in\mst$ and $S\in\mst'_1$ by (\ref{k2}), (\ref{k3}) and Lemma \ref{2.4}. The fact that $\msa_S$ is an intersecting family implies that
\[|\msa_S|=|\msa'_S|\leq \gn{2}{1}=q+1.\] 
Note that $n\geq 2k+4=10$. It follows that
\[\begin{split}
|\msf|&\leq \sum_{S\in\mst'_1}|\msa_S|+\sum_{T\in\mst}|\msf_T\setminus\{D\}|+1\leq |\mst'_1|(q+1)+|\mst|\left(\gn{n-2}{1}-1\right)+1\\
&<(|\mst'_1|+|\mst|)\left(\gn{n-2}{1}-1\right)+1\leq f(n,3,q)
\end{split}\]
by (\ref{eq9}).
$\qed$

\begin{prop}\label{3.8}
Suppose $k\geq 3$ and $n\geq 2k+4$ with $(n,q)\neq (2k+4,2)$. Let $\msf\subseteq \msm(k,n)$ be an intersecting family with $\tau(\msf)=2$. Then $|\msf|\leq f(n,k,q)$. Equality holds iff
\begin{itemize}
	\item[\em(i)] $\msf$ is an HM-type family.
	\item[\em(ii)] $\msf$ is an $\msf_3$-type family when $k=3$.
\end{itemize}
\end{prop}
\proof
It suffices to prove that the upper bounds in Propositions \ref{3.4}--\ref{3.7} are less than $f(n,k,q)$. Let $r=n-2k$. By Lemma \ref{2.10}, it follows that
\[|\msf|/\gn{n-2}{k-2}<1+\frac{q+1}{(q-1)q^{r-1}}\gn{k}{1}\]
in Proposition \ref{3.4}, and
\[|\msf|/\gn{n-2}{k-2}<\left(\frac{1}{q}+\frac{1}{(q-1)q^{r}}\right)\gn{k}{1}+\frac{q^3+q^2-1}{(q-1)q^{r}}\]
in Proposition \ref{3.5} when $m<k$, and
\[|\msf|\leq \gn{3}{1}\gn{n-2}{k-2}\]
in Proposition \ref{3.6} when $k\geq 4$. It is routine to check that these upper bounds are less than
\[\left(1-\frac{1}{(q^2-1)q^r}\right)\gn{k}{1}\gn{n-2}{k-2},\]
which implies the desired result by Lemma \ref{2.11}. 
$\qed$
\begin{prop}\label{3.9}
Suppose $k\geq 3$ and $n\geq 2k+4$ with $(n,q)=(2k+4,2)$. Let $\mathscr{F}\subseteq\mathscr{M}(k,n)$ be a intersecting family with $\tau(\mathscr{F})=t>2$. Then $|\msf|<f(n,k,q)$.
\end{prop}
\proof
Note that $t\leq k$. Let $r=n-2k$. Pick $T\in\mst$. It follows that
\[\msf=\bigcup_{E\in\msm(1,T)}\msf_E,\]
which implies that
\[\begin{split}
|\msf|&\leq\sum_{E\in\msm(1,T)}|\msf_E|\leq\sum_{E\in\msm(1,T)}\gn{k}{1}^{t-1}\gn{n-t}{k-t}\\
&= q^{t-1}\gn{t}{1}\gn{k}{1}^{t-1}\gn{n-t}{k-t}<\frac{q^{t-1}}{((q-1)q^r)^{t-2}}\gn{t}{1}\gn{k}{1}\gn{n-2}{k-2} 
\end{split}\]
by Lemmas \ref{2.3}, \ref{2.14} and \ref{2.10}. Hence, one gets that
\[|\msf|/\gn{k}{1}\gn{n-2}{k-2}<\frac{q^{t-1}(q^t-1)}{(q-1)^{t-1}q^{r(t-2)}}.\]
It is routine to check that
\[\frac{1}{(q^2-1)q^r}+ \frac{q^{t-1}(q^t-1)}{(q-1)^{t-1}q^{r(t-2)}}\leq 1\]
holds for any $t\geq 3$, $q\geq 2$ and $r\geq 4$ with $(q,r)\neq (2,4)$, which implies that $|\msf|<f(n,k,q)$ by Lemma \ref{2.11}.
$\qed$

Together with Proposition \ref{3.8}, we prove Theorem \ref{1.2}.

\end{document}